\documentclass[12pt,a4paper,eqno]{amsart}
\usepackage{graphicx,epsfig}
\usepackage{latexsym}
\usepackage{amssymb}
\usepackage[T1]{fontenc}
\fontencoding{T1}
\usepackage[latin2]{inputenc}
\usepackage{mathtools}
\usepackage{enumerate}
\mathtoolsset{showonlyrefs}
\usepackage{todonotes}

\usepackage{tikz}
\newcommand{\Max}{
  \mathbin{
    \mathchoice
      {\buildMax{\displaystyle}}
      {\buildMax{\textstyle}}
      {\buildMax{\scriptstyle}}
      {\buildMax{\scriptscriptstyle}}
  }
}

\let\oldtocsection=\tocsection
\let\oldtocsubsection=\tocsubsection
\let\oldtocsubsubsection=\tocsubsubsection
\renewcommand{\tocsection}[2]{\hspace{0em}\oldtocsection{#1}{#2}}
\renewcommand{\tocsubsection}[2]{\hspace{1em}\oldtocsubsection{#1}{#2}}
\renewcommand{\tocsubsubsection}[2]{\hspace{2em}\oldtocsubsubsection{#1}{#2}}

\def\Max{\bigcirc\hspace{-4.54mm}\vee \,\,}

\newtheorem{theorem}{Theorem}[section]
\newtheorem{lemma}[theorem]{Lemma}

\newtheorem{definition}[theorem]{Definition}

\newenvironment{namelist}[1]{%
\begin{list}{}
     {
      
      \settowidth{\labelwidth}{#1}
      \setlength{\leftmargin}{1.1\labelwidth}
               }
      }{%
\end{list}}

\newcommand{\eqnsection}{
\renewcommand{\theequation}{\thesection.\arabic{equation}}
 \makeatletter   \csname  @addtoreset\endcsname{equation}{section}
   \makeatother}

\begin{document}
\title[Cramer-Lundberg model for  Markov sequences]
      {Cramer-Lundberg model for some classes of Markov sequences}
\author[Jasiulis-Go{\l}dyn, Lecha\'{n}ska, Misiewicz]
{B.H. Jasiulis-Go{\l}dyn $^1$, A. Lecha\'{n}ska$^{2,*}$ and J.K.~ Misiewicz $^{2,**}$}
\thanks{ $^1$ Institute of Mathematics, University of Wroc{\l}aw, pl. Grunwaldzki 2/4, 50-384 Wroc{\l}aw, Poland, e-mail: Barbara.Jasiulis@math.uni.wroc.pl \\
 $^2$ Faculty of Mathematics and Information Science, Warsaw University of Technology, ul. Koszykowa 75, 00-662 Warszawa, Poland, $^*$ e-mail: lechanska@gmail.com, $^{**}$ e-mail: J.Misiewicz@mini.pw.edu.pl}

\begin{abstract}
The classical  Cramer-Lundberg model  was the first attempt to describe the financial condition of the insurance company. The incomes were approximated by a steady stream of money, insurance payments were not limited and could take any value from zero to infinity. The society did not invest any part of its money, do not have any  employees, shareholders or enterprise maintenance costs. There exists many modifications of the Cramer-Lundberg model which cover at least some of the problems described here, but usually they require  insight into the internal financial policy of the insurance company. We propose here another modification based on Markov processes defined by generalized convolutions. Thanks to the generalized convolutions we can approximate stochastically the internal financial policy of the company based on publicly available data. In this paper we focus on computing the ruin probability for an infinite time horizon for the Markov processes Cramer-Lundberg model where the transition probabilities are defined by generalized convolutions, in particular  $\alpha$-convolution, maximal convolution and the Kendall convolution. 
\end{abstract}
\maketitle


\tableofcontents

\section{Introduction and the classical model description}

\noindent
The classical Cram\'{e}r-L{u}ndberg risk model was introduced by Lundberg in 1903 (see \cite{Lundberg}) and developed by H. Cr\'{a}mer and his Stockholm School at the beginning of the XX century (see e.g. \cite{Cramer}). A reach information on actuarial risk theory, non-life insurance models and  financial models one can find in a huge number of papers and books (e.g. \cite{Franck2, Asmussen, Bowers, Grandel, Klugman, McNeil, Muller, Reiss, Rolski}). We direct the interested reader  to the series of P. Embrechts papers and especially to the book of P. Embrechts, C. Kluppelberg and T. Mikosch \cite{Embrechts}. Between many interesting results we can find there very important modifications of Cramer-Lundberg models for heavy tailed distributions.

The basic model in this theory, called Cramer-Lundberg model or the renewal model, has the following structure:
\begin{namelist}{lll}
\item{\bf (a)} {\bf Claim size process}: the claim sizes $(X_k)$ are i.i.d. positive random variables with cumulative distribution function $F$, finite mean $m = \mathbf{E}X_k$ and variance $\sigma^2 = {\rm Var} X_k$.
\item{\bf (b)} {\bf Claim times}: the claims occur at the random instants of time
    $$
    0 < S_1 < S_2 < \dots \quad a.s.
    $$
    where the inter arrival times
    $$
    T_1 = S_1, \, T_k = S_k - S_{k-1}, \quad k=2,3,\dots
    $$
    are i.i.d. with exponential distribution with mean $\mathbf{E} T_k = {1/{\lambda}}$;
\item{\bf (c)} {\bf Claim arrival process}: the number of claims in the time interval $[0,t]$ is denoted by
    $$
    N(t) = \sup\left\{ n \geqslant 1 \colon S_n < t\right\};
    $$
\item{\bf (d)} the sequences $(X_k)$ and $(T_k)$ are independent of each other.
\end{namelist}
The {\bf risk process} $\{ K_t \colon t \geqslant 0 \}$ is defined by
\begin{equation}
K_t = u + \beta t - \sum_{k=0}^{N_t} X_k,
\end{equation}
where $u \geqslant 0$ denotes the initial capital, $\beta >0$ stands for the premium income rate and $K_t$ is the capital that company have at time $t$.

Notice that this model describes the following situation:
\begin{namelist}{llll}
\item{1.} The insurance company is keeping all the money in a pocket.
\item{2.} All the incomes are coming from the customers insurance premiums
\item{3.} There is no outcomes except for the individual claims of the customers.
\item{4.} There is no cost or benefits coming from the company existence and activities or from the money which the company have.
\end{namelist}

This means that in fact the original simple Cramer-Lundberg model is rather describing
state of water supply in the underground tank when the temperature and humidity in the cave is constant, water drips from the ceiling of the cave with a constant intensity. Animals come to this tank with constant intensity and drink as much water as they need, as long as their needs can be described by some exponential distribution (which is rather reasonable assumption). \\

There are many modifications of the Cramer-Lundberg model. Some of them can be found in \cite{Embrechts}. In these modifications some of the defects  are eliminated. However usually it means that much more information about the insurance company policy is required for description, while no insurance company wants to make such information publicly available. In the next section we propose a family of models based on special class of Markov chains. The secret information of the insurance company can be coded in transition probabilities. This approach may be more convenient than taking under considerations various different company policy elements.\\

The paper is organized as follows: In section 2 we present our model. The considered random walks are very special - the transition probabilities are defined by generalized convolution. A primer on generalized convolutions are described in Section 3, the construction of the random walk with respect to generalized convolution is given in Section 4. The detailed calculations are given in the last 3 sections for the following examples: random walk with respect to stable convolution,  max-convolution and the Kendall random walk. \\

{\bf Notation.} By $\mathbb{N}_0$ we denote here the set of natural numbers including zero if one would have doubts whether zero is a natural number. By  $\mathcal{P}_{+}$ we denote the set of probability measures on the positive half line $[0,\infty)$. If $\lambda_n$ converges weakly to $\lambda$ we write $\lambda_n \rightarrow \lambda$.
For simplicity we will use notation $T_a$ for the rescaling operator (dilatation operator) defined by $(T_a\lambda)(A) = \lambda ({A/a})$ for every Borel set $A$ when $a \neq 0$, and $T_0 \lambda = \delta_0$. \\

\section{Description of the proposed model}
In our model we assume that the insurance company invests at least part of its money, have employees, shareholders which have to have income and at each moment when the insurance payment request comes the company is calculating the total claims amount, subtract from this all costs and add benefits. Thus the corrected cost of the total outcome for claims is not just simple sum of $X_k$. In fact, in this model the financial situation of the company after paying $X_k$ claim can be even better than before. The rich collection of generalized convolutions and freedom in choosing claims distribution $\lambda$ shall give the possibility of adjusting model to the real situation without precise information about company activities. \\[1mm]
We propose here the following structure of the model:
\begin{namelist}{lll}
\item{\bf (a)} {\bf Claim times}: the claims occur at the random instants of time
    $$
    0 <S_1 <S_2 < \dots \quad a.s.
    $$
    where the inter arrival times
    $$
    T_1 = S_1, \, T_k = S_k - S_{k-1}, \quad k=2,3,\dots
    $$
    are iid random variables with exponential distribution, $\mathbf{E} T_k = {1/{\lambda}}$;
\item{\bf (b)} {\bf Claim arrival process}: the number of claims in the time interval $[0,t]$ is the Poisson process with the parameter $\lambda>0$ defined by
    $$
    N(t) = \sup\left\{ n \geqslant 1 \colon S_n < t\right\};
    $$
\item{\bf (c)} {\bf Cumulated claims process}: the total amount of money spent on the first $n$-claims corrected by part of the incomes other than premium and/or some of the costs  is a discrete time Markov process $\{X_n\colon n \in \mathbb{N}_0\}$, which is a $\diamond$-L\'{e}vy process with the step distribution $U_i \sim \mu \in \mathcal{P}_+$ (see Section 3 and 4) and the transition probabilities $P_n(x, \cdot)$. The sequence $(U_i)$ is i.i.d. The cumulative distribution function for the measure $\mu$ we denote by $F$, its density by $f$, the generalized characteristic function by $\Phi_{\mu}$ and we put $H(t) = \Phi_{\mu}(t^{-1})$.
\item{\bf (d)} {\bf Cumulated income units}: the total insurance premium collected by the company up to the moment of $n$-th claim corrected by part of the cost of the company activity and/or part of the income from the investments is a discrete time Markov process $\{ Y_n \colon n \in \mathbb{N}_0\}$,  which is a $\diamond$-L\'{e}vy process with the step distribution $V_i \sim \nu \in \mathcal{P}_+$ (see Section 3 \& 4) and the transition probabilities $\Pi_n(x, \cdot)$. The sequence $(V_i)$ is i.i.d. The cumulative distribution function for the measure $\nu$ we denote by $G$, its density by $g$, the generalized characteristic function by $\Phi_{\nu}$ and we put $J(t) = \Phi_{\nu}(t^{-1})$.
\item{\bf (e)} {\bf Independence assumption:} the processes $\{ N(t) \colon n \in \mathbb{N}_0\}$,  $\{X_n\colon n \in \mathbb{N}_0\}$ and $\{ Y_n \colon n \in \mathbb{N}_0\}$ are independent.
\end{namelist}

The risk process is defined by the following:
$$
R_t = u \oplus \sum_{n=1}^{\infty} Y_n \mathbf{1}_{\{N(t) = n\}}  - \sum_{n=1}^{\infty} X_n \mathbf{1}_{\{N_t = n\}},
$$
where $u \oplus Y_n$ is the Markov process $\{ Y_n \colon n \in \mathbb{N}_0\}$ with the starting point moved to $u>0$ in the generalized convolution sense (see e.g. \cite{BJMR}). Notice that if $N(t)=n$, then we have $R_t = u \oplus  Y_n  -  X_n$, where $X_n = U_1 \oplus U_2 \oplus \cdots \oplus U_n$ and $Y_n = V_1 \oplus V_2 \oplus \cdots \oplus V_n$ and $\oplus$ denotes adding in the generalized convolution sense, i.e $X_n \sim \mu ^{\diamond n}$ and $Y_n \sim \nu ^{\diamond n}$ respectively (for notation $\diamond$ see Section 3).

Let $Q_t(u)$ denotes the probability that the insurance company with the initial capital $u>0$ will bankrupt until time $t$. Since the changes in the process $\{R_t \colon t \geqslant 0\}$ can occur only at the moments $S_n$, $n \geqslant 0$ we see that
\begin{eqnarray*}
Q_t(u) & = &  \mathbf{P}\bigl\{ \exists\, s \leqslant t \colon R_s < 0\bigr\} \\
  & = & 1 - \sum_{n=0}^{\infty} \mathbf{P}\bigl\{ R_{S_k} \geqslant 0 \colon k = 0,1,\dots,n \bigr\} \mathbf{P}\bigl\{N_t = n \bigr\} \\
  & = &  1 - \sum_{n=0}^{\infty} \mathbf{P}\bigl\{u \oplus Y_k >  X_k  \colon k = 0,1,\dots,n \bigr\} \frac{(\lambda t)^n}{n!} e^{-\lambda t}.
\end{eqnarray*}
Calculating the same probability in the unbounded time horizon $Q_{\infty}(u)$ we see that
\begin{eqnarray*}
Q_{\infty} (u) & = &  \mathbf{P}\bigl\{ \exists\, t>0  \colon R_t < 0\bigr\} \\
  & = & 1 - \mathbf{P}\bigl\{u \oplus Y_k >  X_k  \hbox{ for all } k \in \mathbb{N}_0 \bigr\}.
\end{eqnarray*}
Notice that $Q_{\infty}(u)$ does not depend on the process $\{N_t \}$. This is natural since in our case this process is describing only the moments of claims arrival and $\{N_t\}$ is independent of the processes $\{ X_n\}$ and $\{ Y_n\}$. Every continuous time Markov chain taking values in (whole!) $\mathbb{N}_0$ would give the same result.
For abbreviation we introduce the following notation for probability that ruin does not occur:
$$
\delta(u) := Q_{\infty} (u).
$$

\section{Basic information about generalized convolution}
Following K. Urbanik (see \cite{Urbanik64}) we define
\begin{definition}
A commutative and associative $\mathcal{P}$-valued binary
operation $\diamond$ defined on $\mathcal{P}_{+}^2$ is called a
\emph{generalized convolution} if for all $\lambda,\lambda_1,\lambda_2 \in \mathcal{P}_{+}$ and $a \geqslant 0$ we have:
\begin{namelist}{llll}
\item[{\rm (i)}] $\delta_0 \diamond \lambda = \lambda$ ; %
\item[{\rm (ii)}] $(p\lambda_1 +(1-p)\lambda_2) \diamond \lambda = p ( \lambda_1 \diamond \lambda) + (1-p)(\lambda_2
          \diamond \lambda) $ whenever $p\in [0,1]$; %
\item[{\rm (iii)}] $T_a(\lambda_1 \diamond \lambda_2) = (T_a \lambda_1) \diamond (T_a \lambda_2)$ ; %
\item[{\rm (iv)}] if $\lambda_n \rightarrow \lambda$ then $\lambda_n \diamond \eta \rightarrow \lambda \diamond
           \eta$  for all $\eta \in \mathcal{P}$ and $\lambda_n \in \mathcal{P}_{+}$,
\item[{\rm (v)}] there exists a sequence $(c_n)_{n\in\mathbb{N}}$ of positive numbers such that the sequence $T_{c_n}
           \delta_1^{\diamond n}$ converges to a measure different from $\delta_0$.
\end{namelist}
\end{definition}

The set $(\mathcal{P}_{+}, \diamond)$ we call a \emph{generalized convolution algebra}. A continuous mapping $h: \mathcal{P} \to \mathbb{R}$ such that
\begin{itemize}
\item $h(p\lambda + (1-p)\nu)= p h(\lambda) + (1-p) h(\nu)$,
\item $h(\lambda  \diamond \nu)=h(\lambda) h(\nu)$
\end{itemize}
for all $\lambda, \nu \in \mathcal{P}_{+}$ and $p \in (0,1)$, is called a \emph{homomorphism of $(\mathcal{P}_{+}, \diamond)$}. \\

Every convolution algebra $(\mathcal{P}_{+}, \diamond)$ admits two trivial homomorphisms: $h \equiv 1$ and $h \equiv 0$. We say that a generalized convolution is \emph{regular} if it admits a non-trivial homomorphism.  If the generalized convolution is regular then its homomorphism is uniquely determined in the sense that if $h_1, h_2$ are homomorphisms of $(\mathcal{P}_{+}, \diamond)$ then there exists $c>0$ such that $h_1(\lambda) = h_2(T_c \lambda)$ (for details see \cite{Urbanik64}). It was also shown in \cite{Urbanik64} that the generalized convolution is regular if and only if there exists unique up to a scale function
$$
\mathcal{P}_{+} \ni \lambda \longrightarrow \Phi_{\lambda} \in C([0,\infty))
$$
such that  for all $\lambda, \nu, \lambda_n \in \mathcal{P}_{+}$ the following conditions hold:
\begin{namelist}{ll}
\item{1.} $\Phi_{p \lambda + q \nu}(t) = p \Phi_{\lambda}(t) + q \Phi_{\nu}(t)$, for $p,q \geqslant 0$, $p+q = 1$;
\item{2.} $\Phi_{\lambda \diamond \nu}(t) = \Phi_{\lambda}(t) \Phi_{\nu}(t)$;
\item{3.} $\Phi_{T_a \lambda}(t) = \Phi_{\lambda}(at)$ for $a \geqslant 0$;
\item{4.} the uniform convergence of $\Phi_{\lambda_n}$ on every compact set to a function $\Phi$ is equivalent with the existence of $\lambda \in \mathcal{P}_+$ such that $\Phi = \Phi_{\lambda}$ and $\lambda_n \rightarrow \lambda$
\end{namelist}
The function $\Phi_{\lambda}$ is called the $\diamond$-generalized characteristic function of the measure $\lambda$. Let $\Omega(t) = h(\delta_t)$. By properties 1 and 2 of the  characteristic function we see that
$$
\Phi_{\lambda}(t) = \int_0^{\infty} \Omega(xt) \lambda(dx),
$$
thus the function $\Omega$ is called the kernel of generalized characteristic function (similarly as the function $e^{it}$ is the kernel of Fourier transform, i.e. the classical characteristic function).

\vspace{2mm}

\noindent
{\bf Examples.} For details see \cite{Bingham, CKS, JasKula, King, KU2, JM_l-c, Urbanik64, Urbanik76,Urbanik88}.
\begin{namelist}{ll}
\item[{\bf 3.0.}] The {\em classical convolution}, denoted by $\ast$ is given by:
    $$
    \delta_a \ast \delta_b = \delta_{a+b}.
    $$
    Here we have $\Omega(t) = e^{-t}$ if we consider this convolution on $\mathcal{P}_{+}$ and $\Omega(t) = e^{it}$
    if we consider it on the whole line.
\item[{\bf 3.1.}] {\em Symmetric convolution} on $\mathcal{P}_{+}$ is defined by
    $$
    \delta_a \ast_{s} \delta_b = \frac{1}{2}\, \delta_{|a-b|} + \frac{1}{2} \,\delta_{a+b}.
    $$
    The kernel of generalized characteristic function here is
    $\Omega(t) = \cos(t)$.
\item[{\bf 3.2.}] By {\em stable convolution} $\ast_{\alpha}$ for  $\alpha>0$ we understand the following:
    $$
    \delta_a \ast_{\alpha} \delta_b = \delta_c, \quad c = (a^{\alpha}+b^{\alpha})^{1/{\alpha}}, \quad a,b \geqslant 0.
    $$
    The kernel of generalized $\ast_{\alpha}$-characteristic function is $\Omega(t) = e^{-t^{\alpha}}$.
\item[{\bf 3.3.}] {\em $\infty$-convolution} is defined by
$$
\delta_a \Max \delta_b = \delta_{\max\{a,b\}}.
$$
This convolution admits existence of characteristic function, but its kernel is not  continuous: $\Omega(t) = \mathbf{1}_{[0,1]}(t)$.
\item[{\bf 3.4.}] The {\em Kendall convolution} $\vartriangle_{\alpha}$ on $\mathcal{P}_{+}$, $\alpha > 0$, is defined by
$$
\delta_x \vartriangle_{\alpha} \delta_1 = x^{\alpha} \pi_{2\alpha} + (1-x^{\alpha}) \delta_1, \quad x\in [0,1],
$$
where $\pi_{2\alpha}$ is a Pareto measure with density $2\alpha x^{-2\alpha -1} \mathbf{1}_{[1,\infty)}(x)$.
The kernel of generalized characteristic function  here is given by $\Omega(t) = \left( 1 - (ts)^{\alpha}\right)_{+}$,
   where $a_{+} = a$ for $a\geqslant 0$ and $a_{+} = 0$ otherwise.
\item[{\bf 3.5.}] The {\em Kingman convolution} $\otimes_{\omega_s}$ on $\mathcal{P}_{+}$, $s> - \frac{1}{2}$, is defined by
$$
\delta_a \otimes_{\omega_s} \delta_b = \mathcal{L}\left(\sqrt{ a^2 + b^2 + 2ab \theta_s }\right),
$$
where $\theta_s$ is absolutely continuous with the density function
$$
f_s (x)= \frac{\Gamma(s+1)}{\sqrt{\pi}\, \Gamma(s + \frac{1}{2})} \bigl( 1 - x^2\bigr)_{+}^{s - \frac{1}{2}}.
$$
The kernel of generalized characteristic function here is given by the Bessel function of the first kind with parameter connected with $s$.

\item[{\bf 3.6.}]
For every $p \geqslant 2$ and properly chosen $c>0$  the function $h(\delta_t) = \varphi(t) = \varphi_{c,p}(t) = ( 1 - (c+1)t +ct^p)\mathbf{1}_{[0,1]}(t)$ is the kernel of a Kendall type (see \cite{JM_l-c}) generalized  convolution $\diamond$ defined for $x \in [0,1]$ by the formula:
$$
\delta_x \diamond \delta_1 = \varphi(x) \delta_1 + x^p \lambda_1 + (c+1)(x-x^p) \lambda_2,
$$
where $\lambda_1, \lambda_2$ are probability measures absolutely continuous with respect to the Lebesgue measure and independent of $x$. For example if $c = (p-1)^{-1}$ then
$$
\lambda_1(du) = \frac{2c}{u^3} \Bigl[ (c+1)(p+1) u^{1-p} + (c+1)(p-2) + cp (2p-1) u^{-2p-2} \Bigr] \!\!\mathbf{1}_{[1,\infty)}(u) du,
$$
and
$$
\lambda_2(du) = c \bigl[ 2(p-2) + (p+1)u^{-p+1} \bigr] u^{-3} \mathbf{1}_{[1,\infty)}(u) du.
$$
\end{namelist}

\section{Random walk with respect to the generalized convolution}
All the information contained in this section comes from \cite{BJMR}, where the L\`{e}vy processes with respect to generalized convolution were defined and studied. It was shown there that each such process is a Markov process (in the classical sense) with the transition probabilities defined by generalized convolution. We consider here only discrete time stochastic processes of this kind.

\begin{definition} A discrete time stochastic process $\left\{ X_n \colon n \in \mathbb{N}_0 \right\}$ is a random walk with respect to generalized convolution $\diamond$ with the step distribution $\mu$ if it is the Markov process with the transition probabilities
$$
P_{k,n} (x, dy) = \delta_x \diamond \mu^{\diamond (n-k)}(dy), \quad n \geqslant k.
$$
\end{definition}

The consistency of this definition and the existence of the random walk with respect to generalized convolution $\diamond$ was shown in \cite{BJMR}. Notice that in the case of classical convolution it is the simple random walk with the step distribution $\mu$ and it can be simply represented as $X_n = U_1 + \dots + U_n$, where $(U_k)$ is a sequence of i.i.d. random variables with distribution $\mu$.

There are only two cases, when generalized convolution $\diamond$ is representative, i.e. there exists a sequence of functions $f_n \colon R^n \rightarrow R$ such that $X_n = f_n(U_1, \dots, U_n)$: \\
{\bf 4.2.} for the $\ast_{\alpha}$-convolution  $X_n = \bigl( U_1^{\alpha} + \dots + U_n^{\alpha} \bigr)^{1/{\alpha}}$, \\
{\bf 4.3.} for the $\infty$-convolution we have $X_n = \max\{ U_1, \dots, U_n\}$. \\
For other generalized convolutions rewriting convolution in the language of the corresponding independent random variables is more complicated (if possible) and requires assistance of some additional variables. For example we have \\
{\bf 4.4.} for the Kendall convolution for $x \in [0,1]$ the measure $\delta_x \diamond_{\alpha} \delta_1$ is the distribution of the random variable
$$
\bigl(x \oplus_{\vartriangle_{\alpha}} 1 \bigr)(\omega) := \mathbf{1}_{\{Q(\omega) > x^{\alpha}\}} + \mathbf{1}_{\{Q(\omega) \leqslant x^{\alpha}\}} \Pi_{2\alpha} (\omega),
$$
where $Q$ has uniform distribution on $[0,1]$, $\Pi_{2\alpha}$ has the Pareto distribution with the density $\pi_{2\alpha}$ described in example 3.4, $Q$ and $\Pi_{2\alpha}$ are independent; \\
{\bf 4.5.} for the Kingman convolution and $a,b>0$ we can define:
$$
\bigl(a \oplus_{\omega_s} b\bigr)(\omega) := \sqrt{a^2 + b^2 + 2ab \theta_s },
$$
where $\theta_s$ is absolutely continuous with the density function $f_s$ described in example 3.5.

\section{Model for $\ast_{\alpha}$ random walk}

For $\ast_{\alpha}$ generalized convolution on $\mathcal{P}_{+}$ we have
$$
X_n = \left( X_{n-1}^{\alpha} + U_n^{\alpha} \right)^{1/{\alpha}} = \left( U_1^{\alpha} + \dots + U_n^{\alpha} \right)^{1/{\alpha}}, \quad n \geqslant 1,
$$
where $(U_k)$ are independent identically distributed random variables with cumulative distribution function $F_U$ responsible for the damage claim values. By $F = F_{U^\alpha}$ we denote the cumulative distribution function of $U^\alpha$. We assume also that $m_{\alpha} = \mathbf{E} U_1^{\alpha} < \infty$.

We assume here that the variables $V_k$, responsible for the insurance premium during the time $T_k$ are independent identically distributed with the cumulative distribution function  $F_{V}(x) = 1 - e^{-\gamma x^{\alpha}}$, for $x>0$. This assumption seems to be natural, since this is the distribution with the lack of memory property (see \cite{Poisson}) for $\ast_{\alpha}$-convolution.  Consequently
$$
Y_n = \left( Y_{n-1}^{\alpha} + V_n^{\alpha} \right)^{1/{\alpha}} = \left( V_1^{\alpha} + \dots + V_n^{\alpha} \right)^{1/{\alpha}}, \quad n \geqslant 1.
$$
Now we have
$$
\mathcal{R}_t = \bigg[ u^\alpha+\beta^\alpha \sum\limits_{n=1}^{\infty} Y_n^\alpha\mathbf{1}_{\{N_{t} = n\}} \bigg]^{\frac{1}{\alpha}}
 -  \bigg[\sum\limits_{n=1}^{\infty} X_n^\alpha \mathbf{1}_{\{N_t = n\}} \bigg]^{\frac{1}{\alpha}}.
$$

We want to calculate the ruin probability (see \cite{Asmussen, Embrechts}) for the insurance company  $u$ by the time $t$:
$$
Q_{t}(u) = P \bigl\{ \exists s \leq t: \mathcal{R}_s \leqslant 0 \bigr\}
$$
in the special case $t = \infty$, i.e.:
$$
Q_{\infty}(u)= P\bigl\{ \exists t > 0: \mathcal{R}_t \leqslant 0 \bigr\}.
$$
Since the ruin can occur only in the claims arrival moments i.e. in the moments of jumps of the Poisson process ${N}_t$, thus it is enough to consider $\mathcal{R}_{S_n}$:
\begin{eqnarray*}
Q_{\infty}(u) & = &  \lefteqn{P\bigl\{ \mathcal{R}_t = 0 \;\; \text{ for some }  t>0 \bigr\}
= 1 - P \bigl\{ \mathcal{R}_{S_n}>0 \;\;\;  \forall n \in \mathbb{N} \bigr\}} \\
& = & 1 - P \biggl\{ u^\alpha - \sum\limits_{i=1}^n \bigl( U_i^\alpha - \beta^\alpha V_i ^\alpha \bigr) > 0 \;\;\; \forall n\in\mathbb{N} \biggr\}\\
& = & 1 - P \biggl\{ \sup_{n \geqslant 1} \sum\limits_{i=1}^n \bigl( U_i^\alpha-\beta^\alpha V_i^\alpha \bigr) < u^\alpha \biggr\} =: 1 - \delta(u^\alpha).
\end{eqnarray*}
Basically, the function $\delta(u^\alpha)$ we can calculate following the classical calculations:
\begin{eqnarray*}
\lefteqn{\delta ( u^\alpha ) = \mathbf{P} \Big\{ \sup_{n \geqslant 2} \sum\limits_{i=2}^n (U_i^\alpha-\beta^\alpha V_i^\alpha ) < u^\alpha + \beta^\alpha V_1^\alpha - U_1^\alpha, \,\, U_1^\alpha - \beta^\alpha V_1^\alpha < u^\alpha \Big\}} \\
&=&\int\limits_0^\infty  \gamma e^{-\gamma y} \! \int\limits_0^{u^\alpha + \beta^\alpha y}\! P \Bigl\{ \sup_n  \sum\limits_{i=2}^n ( U_i^\alpha - \beta^\alpha V_i^\alpha ) < u^\alpha + \beta^\alpha y - x \Bigr\} \, F_x(\mathrm{d}x) \mathrm{d}y\\
&=& \beta^{-\alpha } e^{ \frac{\gamma u^\alpha}{\beta^\alpha} }  \int\limits_{u^\alpha}^\infty\! \gamma e^{\frac{-\gamma z}{\beta^\alpha}} \int\limits_0^z\!\delta(z-x)\,\mathrm{d}F(x) \mathrm{d}z.
\end{eqnarray*}
In the last step in these calculations we substituted $u^\alpha+\beta^\alpha y=z$.
Now we calculate the derivative of both sides of this equality with respect to  $\mathrm{d} u^\alpha $:
$$
\frac{\mathrm{d}\delta(u^\alpha)}{\mathrm{d}u^\alpha}=\frac{\gamma}{\beta^\alpha} \, \delta(u^\alpha)-\frac{\gamma}{\beta^\alpha}\int\limits_0^{u^\alpha}\delta(u^\alpha-x)\mathrm{d}F(x).
$$
Integrating both sides of this equality over the set  [0, t] with respect to the measure with the density function  $\alpha u^{\alpha-1}$ for $u>0$ we obtain:
\begin{align*} \delta(t^\alpha)=\delta(0)+\frac{\gamma\alpha}{\beta^\alpha}\int\limits_0^t u^{\alpha-1}\delta(u^\alpha)\mathrm{d}u
-\frac{\gamma\alpha }{\beta^{\alpha}}\int\limits_0^t\int\limits_0^{u^\alpha} u^{\alpha-1}\delta(u^\alpha-x)\mathrm{d}F(x)\mathrm{d}u.
\end{align*}
The first integral on the right hand side we denote by $I_1$, second by $I_2$. Then
$$
I_1 = \alpha \int\limits_0^{t} u^{\alpha-1} \delta(u^{\alpha} )\mathrm{d}u = \int\limits_0^{t^\alpha} \delta(y)\mathrm{d}y = \int\limits_0^{t^\alpha} \delta(t^\alpha - x)\mathrm{d}x.
$$
In the second integral we change order of integration and then substitute $u^\alpha - x = r$:
\begin{eqnarray*}
\lefteqn{I_2 = \int\limits_0^t\int\limits_0^{u^\alpha} \alpha u^{\alpha-1} \delta(u^\alpha-x) \mathrm{d}F(x) \mathrm{d}u
= \int\limits_0^t \int\limits_{x^{\frac{1}{\alpha}}}^t \alpha u^{\alpha-1}\delta(u^\alpha-x) \mathrm{d}u \, \mathrm{d}F(x) } \\
& = & \int\limits_0^{t^\alpha} \int\limits_0^{t^\alpha-x} \! \delta(r)\mathrm{d} \, \mathrm{d}F(x)
 = \biggl[ F(x) \!\! \int\limits_0^{t^\alpha-x} \! \delta(r)\mathrm{d}r   \biggr]_{x=0}^{t^\alpha}+\int\limits_0^{t^\alpha} \! \delta(t^\alpha-x)F(x)\mathrm{d}x.
\end{eqnarray*}
The last equality we obtained integrating by parts. Since F(0)=0 we obtain
$$
\delta(t^\alpha)=\delta(0)+\frac{\gamma}{\beta^\alpha}\int\limits_0^{t^{\alpha}} \delta(t^\alpha-x)\left(1-F(x)\right)\mathrm{d}x.
$$
In order to calculate $\delta(0)$ notice first that ruin probability for the insurance company with infinite initial capital is zero, thus $\delta(\infty) = 1$ and we have
$$
1=\delta(0)+\frac{\gamma}{\beta^\alpha}\int\limits_0^\infty(1-F(x))\mathrm{d}x= \delta(0)+ \frac{\gamma}{\beta^\alpha}  \mu_{\alpha}.
$$
Thus $\delta(0)= 1 - \frac{\gamma}{\beta^\alpha}\mu_{\alpha}$ and we have
$$
\delta(t^\alpha)=1-\frac{\gamma}{\beta^\alpha} \mu_{\alpha} + \frac{\gamma}{\beta^\alpha} \int\limits_0^{t^\alpha}\delta(t^\alpha-x)(1-F(x))\mathrm{d}x.
$$
For the convenience in further calculations we substitute $t^{\alpha}= z$. Let $\widehat{f}$ be  the Laplace-Stietjes transform given by $\widehat{f}(s) = \int\limits_0^\infty e^{-zs} f(z)\mathrm{d}z$. Thus
\begin{eqnarray*}
\lefteqn{\hspace{-20mm} \int\limits_0^\infty e^{-zs} \int\limits_0^z\delta(z-x)G(x)\mathrm{d}x\mathrm{d}z
= \int\limits_0^\infty \int\limits_x^\infty e^{-zs}\delta(z-x)\,\mathrm{d}z\, G(x)\,\mathrm{d}x } \\
 & = &\int\limits_0^\infty\int\limits_0^\infty e^{-(x+z)s}\delta(z)\, \mathrm{d}z\, G(x)\, \mathrm{d}x
 = \widehat{\delta}(s)\widehat{G}(s).
\end{eqnarray*}
Consequently we obtain:
$$
\widehat{\delta}(s)=\Bigl(1-\frac{\gamma}{\beta^\alpha} \mu_{\alpha} \Bigr) \frac{1}{s}+\frac{\gamma}{\beta^\alpha}\widehat{\delta}(s)\widehat{G}(s),
$$
thus
$$
\widehat{\delta}(s) =\frac{\beta^\alpha  -\gamma \mu_{\alpha}}{(\beta^\alpha - \gamma \widehat{G}(s))s}. \eqno{(\ast)}
$$
Since the Laplace-Stietjes transform uniquely determines function, we finally  have that in this case the ruin probability for the insurance company with the initial capital $u$ is equal $Q_{\infty}(u) = 1 - \delta(u^{\alpha})$ with the function $\delta$ obtained from the equation ($\ast$).

\section{Model for $\infty$-generalized convolution}

For the random walk with respect to the $\infty$-convolution on $\mathcal{P}_{+}$ we have $X_n = \max\{U_1, \dots, U_n \}$ and $Y_n = \max\{V_1, \dots, V_n \}$, where $(U_k)$ and $(V_k)$ are independent sequences of i.i.d. positive random variables with distributions $\mu$ and $\nu$ and the cumulative distribution functions $F$ and $G$ respectively.
Consequently $X_n$ has the cumulative distribution function $F^n$, $Y_n$ has this function equal $G^n$ and $u \oplus Y_n$ has $G^n$.

The first safety condition for the insurance company is $\mathbf{E} R_t > 0$, thus we need to calculate $\mathbf{E} X_t$ and $\mathbf{E} (u \oplus Y_t)$, where
$$
X_t = \sum_{n=1}^{\infty} X_n \mathbf{1}_{\{N_t = n\}}, \quad Y_t = \sum_{n=1}^{\infty} Y_n \mathbf{1}_{\{N(t) = n\}}.
$$
We have
\begin{eqnarray*}
\mathbf{E} X_t & = & \sum_{n=0}^{\infty} \int_0^{\infty} \!\! x \;\frac{(\lambda t)^{n}}{n!} e^{-\lambda t}\,  d_x \left[ F(x)^n \right] =  e^{-\lambda t} \int_0^{\infty}\!\!  x \, d_x \left[ \exp\{ \lambda t F(x)\} \right] \\
& = & \lambda t \int_0^{\infty}  x\, e^{-\lambda t \overline{F}(x) } dF(x),
\end{eqnarray*}
where $\overline{F} = 1 - F$ is the survival function for $U_k$. In order to calculate $\mathbf{E} (u \oplus Y_t)$ notice first that the variable $u \oplus Y_n$ is taking value $u$ with probability $G(u)^n$, thus
$$
\mathbf{E} (u \oplus Y_t) = \lambda t u e^{- \lambda t \overline{G}(u)} + \lambda t \int_{(u, \infty)}\!\! x\, e^{-\lambda t \overline{G}(x) } dG(x)
$$
with the same notation $\overline{G} = 1 - G$. Consequently the first safety condition in the case of $\infty$-generalized convolution is the following:
$$
u e^{- \lambda t \overline{G}(u)} + \int_{(u, \infty)}\!\! x\, e^{-\lambda t \overline{G}(x) } dG(x) \geqslant \int_0^{\infty}  x\, e^{-\lambda t \overline{F}(x) } dF(x).
$$
Usually we take the random variables $V_k$ with the distribution having the lack of memory property, which in the case of $\infty$-convolution is given by the cumulative distribution function $G(x) = \mathbf{1}_{(a,\infty)}(x)$ for some $a>0$ (see \cite{Poisson} for details). In this case we have
$$
\mathbf{E} (u \oplus Y_t) = \lambda t \bigl( u\vee a \bigr)  e^{- \lambda t \overline{G}(u\vee a)}
$$
Consequently the first safety condition for the $\infty$-convolution is the following:
$$
\bigl( u\vee a \bigr)  e^{- \lambda t \overline{G}(u\vee a)} - \int_0^{\infty}  x\, e^{-\lambda t \overline{F}(x) } dF(x) \geqslant 0 \quad \quad \forall \, t> 0.
$$
Calculating the probability that the company will not bankrupt in the unbounded time horizon we shall consider two cases. If the insurance premium has distribution with the lack of memory property, i.e. $G(x) = \mathbf{1}_{(a,\infty)}(x)$, then we have $Y_k = V_1=a$ for all $k \in \mathbb{N}$, thus
\begin{eqnarray*}
\delta(u) & = & \mathbf{P} \bigl\{ u \oplus Y_k > X_k \hbox{ for all } k \in \mathbb{N} \bigr\} = \mathbf{P} \bigl\{ u \vee a > X_k \hbox{ for all } k \in \mathbb{N} \bigr\}\\
  & = & \mathbf{P} \left( \bigcap_{k\in \mathbb{N}} \left\{ X_1 < u \vee a, \dots, X_k < u \vee a  \right\}\right) \\
  & = & \lim_{k\rightarrow \infty} \mathbf{P} \bigl\{ X_1 < u \vee a, \dots, X_k < u \vee a  \bigr\} = \lim_{k\rightarrow \infty} F^k \bigl(u \vee a \bigr).
\end{eqnarray*}
We see that bankruptcy in unbounded time horizon is granted if only the random variables $U_k$ can take any positive value, i.e. if $F(x) <1$ for all $x>0$. However if the biggest possible claim is less than $u \vee a$ then bankruptcy is impossible and $\delta(u) = 1$.

If we assume that the cumulative distribution functions $F, G$ are not trivial  then we have
\begin{eqnarray*}
\delta(u) & = & \mathbf{P} \left\{ u \oplus Y_k > X_k\, \hbox{ for all }\, k \in \mathbb{N} \right\} \\
  & = & \mathbf{P} \left\{ u \vee V_1 > X_1,\,  u \vee V_1 \vee Y_k' > U_1 \vee X_k'\, \hbox{ for all } k \in \mathbb{N} \right\} \\
  & = & \int_{\{u \vee y > x\}} \mathbf{P} \left\{u \vee y \vee Y_k' >x \vee X_k'\, \hbox{ for all } k \in \mathbb{N} \right\} \, dF(x) \, dG(y) \\
  & = & \int_{\{u \vee y > x\}} \mathbf{P} \left\{u \vee y \vee Y_k' > X_k'\, \hbox{ for all } k \in \mathbb{N} \right\} \, dF(x) \, dG(y) \\
  & = & \int_0^u \int_0^u \delta(u) \, dF(x) \, dG(y) + \int_u^{\infty} \int_0^y \delta(y) \, dF(x) \, dG(y),
\end{eqnarray*}
where $Y_k' = \max\{ V_2, \dots, V_{k+1} \}$ and $X_k' = \max\{ U_2, \dots, U_{k+1}\}$.
For $a>b>0$ let $\delta(a,b) = \mathbf{P} \left\{a \vee Y_k' > b \vee X_k'\, \hbox{ for all } k \in \mathbb{N} \right\}$. Thus we can write:
$$
\delta(u) = \delta(u) G(u) F(u) + \int_u^{\infty} \delta(y) F(y) \, dG(y).
$$
If the distribution functions $F,G$ have densities $f,\,g$ then differentiating both sides of the previous equation we obtain
$$
\delta'(u) = \delta'(u) F(u) G(u) + \delta(u) G(u) f(u).
$$

{\bf Example.} Assume that for $0<a<b$ we have
$$
F(x) = \frac{x}{a} \; \mathbf{1}_{(0,a]}(x) + \mathbf{1}_{(a,\infty)}(x), \quad G(x) = \frac{x}{b}\; \mathbf{1}_{(0,b]}(x) + \mathbf{1}_{(b,\infty)}(x).
$$
Then for $u \in (0,a)$
$$
\frac{\delta'(u)}{\delta(u)} = \frac{\frac{u}{ab}}{ 1 - \frac{u^2}{ab}}, \quad \hbox{thus } \quad \delta(u) = \frac{\delta(0)}{ \sqrt{1 - \frac{u^2}{ab}}}.
$$
For $u \in (a,b)$ we have
$$
\delta'(u) = \delta'(u) \frac{x}{b}\; \mathbf{1}_{(a,b]}(u)\quad \hbox{thus}\quad \delta(u) = \hbox{const}.
$$
If $u>b$ then evidently $\delta(u) = 1$. Since in our case the function $\delta$ is continuous then $\delta(0) = \sqrt{1 - {a/b}}$, thus finally
$$
\delta(u) = \sqrt{ \frac{1 - \frac{a}{b}}{ 1 - \frac{u^2}{ab}}}.
$$
\section{Model for the Kendall random walk}

As we have seen in the previous section, the renewal process based on the $\max$-generalized convolution is rather trivial, for example in the case of very natural step distribution $F(x) = \mathbf{1}_{(1,\infty)}(x)$ it is not moving at all. The Kendall convolution generalizes the $\max$-convolution in the sense that the Kendall convolution of non-negative two random variables is equal to the maximum one with positive probability, otherwise it is is larger than the maximum. For this generalized convolution we will not get any trivial process. \\

The Kendall random walk, i.e. random walk with respect to Kendall convolution,  $\{ X_n\colon n \in \mathbb{N}_0 \}$ for fixed $\alpha>0$ can be described by the recursive construction given below. We see that we can get here explicit formulas for $X_n$, but except the sequence $(U_n)$ we need also two sequences of random variables (catalyzers of $\vartriangle_{\alpha}$-adding)
\begin{namelist}{lll}
\item{1.} $(U_k)$ i.i.d. random variables with distribution $\mu$;
\item{2.} $(\xi_k)$ i.i.d. random variables with uniform distribution on $[0,1]$;
\item{3.} $(\Pi_k)$ i.i.d. random variables with distribution Pareto $\pi_{2\alpha}$;
\end{namelist}
where all these sequences are independent. Then the Kendall random walk has the following representation: $X_0 \equiv 0$,
$$
X_1 = U_1, \quad \quad X_{n+1} = M_{n+1}\bigl[ \mathbf{1}(\xi_{n+1} > \varrho_{n+1}) + \Pi_{n+1} \mathbf{1}(\xi_{n+1} < \varrho_{n+1})\bigr],
$$
where
$$
M_{n+1} = \max\{ X_n, U_{n+1}\}, \quad m_{n+1} = \min\{ X_n, U_{n+1}\}, \quad \varrho_{n+1} = \frac{m_{n+1}^{\alpha}}{M_{n+1}^{\alpha}}.
$$
This representation is especially helpful if we want to make computer simulation of the Kendall random walk. For calculations however it is more convenient to use the Markov properties and transition probabilities given in Definition 4.1.\\

We consider here two Markov chains $\{ X_n\colon n \in \mathbb{N}_0 \}$ and $\{ Y_n\colon n \in \mathbb{N}_0 \}$ with transition probabilities given respectively as follows:

\begin{lemma}\label{lem:1}
For all $x,y,t \geq 0$ and $\mu,\nu \in \mathcal{P}_+$ we have
\begin{eqnarray*}
&& \hspace{-15mm} h(x,y,t) := \delta_x \vartriangle_{\alpha} \delta_y (0,t) = \left( 1 - \left( \frac{xy}{t^2}\right)^{\alpha} \right) \mathbf{1}_{\{x < t, y < t\}} \\
 & &  = \left[ \Psi\left(\frac{x}{t}\right) +  \Psi\left(\frac{y}{t}\right) - \Psi\left(\frac{x}{t}\right) \Psi\left(\frac{y}{t}\right)\right] \mathbf{1}_{\{x < t, y < t\}}, \\
 & &  \hspace{-15mm} \left( \delta_v \vartriangle_{\alpha} \mu \right) \,(0,t)= P_1 (v, (0,t)) \\
 & &  = \left[\Psi \left( \frac{v}{t} \right) F(t) + \left(1-\Psi \left( \frac{v}{t} \right)\right) H(t)\right] \mathbf{1}_{\{v < t \}},\\
&&  \hspace{-15mm}  \left( \delta_u \vartriangle_{\alpha} \nu \right) \,(0,t) = \Pi_1 (u, (0,t)) \\
 & & = \left[\Psi \left( \frac{u}{t} \right) G(t) + \left(1-\Psi \left( \frac{u}{t} \right)\right) J(t)\right] \mathbf{1}_{\{u < t \}},
\end{eqnarray*}
where $F(t) = \mu \, (0,t], G(t) = \nu \, (0,t]$ and
$$
H(t) = \int\limits_0^t \Psi \left( \frac{y}{t}\right) \mu(dy), \quad 
J(t) = \int\limits_0^t \Psi \left( \frac{y}{t}\right) \nu(dy).
$$
\end{lemma}

\noindent{\bf Proof.} 
For the second formula of the first equality it is sufficient to apply $\Psi\left(\frac{x}{t}\right)=\left(1-\frac{x^{\alpha}}{t^{\alpha}}\right)_+$.
The second equation one can prove for $v < t$ in the following way:
\begin{eqnarray*}
 \left( \delta_v \vartriangle_{\alpha} \mu \right) \,(0,t) & = &\int\limits_0^{t} \left[ \Psi\left(\frac{v}{t}\right) + \Psi\left(\frac{y}{t}\right) - \Psi\left(\frac{v}{t}\right) \Psi\left(\frac{y}{t}\right)\right]  \mu(dy) \\
 & = & \left[\Psi \left( \frac{v}{t} \right) F(t) + \left(1-\Psi \left( \frac{v}{t} \right)\right) H(t)\right] \mathbf{1}_{\{v < t \}}.
 \end{eqnarray*}
 The last formula is equivalent with the previous one.
\qed

\begin{lemma}\label{lem:2}
Let  $\mu, \nu\in \mathcal{P}_+$ and $n \in \mathbb{N}$. Then $\left(\delta_u \vartriangle_{\alpha} \nu^{\vartriangle_{\alpha} n}\right) \,(0,t) = \left(\delta_v \vartriangle_{\alpha} \mu^{\vartriangle_{\alpha} n}\right) \,(0,t) = 0$ for $u,v >t$ and for $u,v < t$ we have 
\begin{eqnarray*}
 \left(\delta_v \vartriangle_{\alpha} \mu^{\vartriangle_{\alpha} n}\right) \,(0,t) & = &  F_{v,n}(t) = \left[\Psi \left( \frac{v}{t} \right) F_n(t) + \left(1-\Psi \left( \frac{v}{t} \right)\right) H_n(t) \right] \\
& = & H(t)^{n-1} \left[n \left(F(t)-H(t)\right) \Psi \left( \frac{v}{t} \right) + H(t)\right]
\end{eqnarray*} 
and
\begin{eqnarray*}
 \left(\delta_u \vartriangle_{\alpha} \nu^{\vartriangle_{\alpha} n}\right) \,(0,t) & = &  G_{u,n}(t) = \left[\Psi \left( \frac{u}{t} \right) G_n(t) + \left(1-\Psi \left( \frac{u}{t} \right)\right) J_n(t) \right] \\
& = & J(t)^{n-1} \left[n \left(G(t)-J(t)\right) \Psi \left( \frac{u}{t} \right) + J(t)\right],
\end{eqnarray*}
where $F_n(t) = \mu^{\vartriangle_{\alpha} n} (0,t], G_n(t) = \nu^{\vartriangle_{\alpha} n}(0,t]$ (with notation
$F:= F_0, \; G:= G_0,\; H:= H_0, J:= J_0$) and
$$
H_n(t) = \int\limits_0^t \Psi \left( \frac{y}{t}\right) \mu^{\vartriangle_{\alpha} n} (dy), \quad 
J_n(t) = \int\limits_0^t \Psi \left( \frac{y}{t}\right) \nu^{\vartriangle_{\alpha} n} (dy).
$$

\end{lemma}

\noindent
{\bf Proof.}
By the second formula from Lemma \ref{lem:1} we have
$$
\left(\delta_v \vartriangle_{\alpha} \mu^{\vartriangle_{\alpha} n}\right) \,(0,t) = \left[\Psi \left( \frac{v}{t} \right) F_n(t) + \left(1-\Psi \left( \frac{v}{t} \right)\right) H_n(t)\right] \mathbf{1}_{\{v < t \}}.
$$
 Since
 $$
 H_n(t) = \int\limits_0^t \Psi \left( \frac{y}{t}\right) \mu^{\vartriangle_{\alpha} n} (dy) = 
 \left(\int\limits_0^t \Psi \left( \frac{y}{t}\right) \mu (dy)\right)^n = H(t)^n
 $$
 and
 $$
F_n(t) =   H(t)^{n-1} \bigl[ H(t) + n(F(t) - H(t))\bigr],
$$
we arrive to 
$$
 \left(\delta_v \vartriangle_{\alpha} \mu^{\vartriangle_{\alpha} n}\right) \,(0,t) = H(t)^{n-1} \left[n \left(F(t)-H(t)\right) \Psi \left( \frac{v}{t} \right) + H(t)\right],
$$
which ends the proof. \qed

\subsection{Inversion formula and cumulative distribution functions} 

The generalized characteristic function for the Kendall convolution is the William\-son integral transform:
\begin{eqnarray}
\Phi_{\mu}(t) & :\stackrel{(1)}{=} & \int_0^{\infty} \left( 1 - (ts)^{\alpha} \right)_{+} d F(s) \\
& \stackrel{(2)}{=} & F(1/t) - t^{\alpha} \int_0^{1/t} s^{\alpha} dF(s) \stackrel{(3)}{=}  \alpha t^{\alpha} \int_0^{1/t} s^{\alpha - 1} F(s) ds.
\end{eqnarray}
Notice that the Williamson transform (see \cite{AJO,KendallWalk,factor,renewal, Williamson}) is easy to invert:

If $\mu$ has cumulative distribution function $F$, then for $H(t) := \Phi_{\mu}(1/t)$ using the formulation (3) we have
$$
 \int_0^{t} s^{\alpha - 1} F(s) ds = \alpha^{-1} t^{\alpha} H(t).
$$
Differentiating both sides with respect to $t$ we obtain
$$
F(t) = \alpha^{-1} t^{1- \alpha} \frac{d}{dt} \left( t^{\alpha} H(t)\right) = H(t) + \alpha^{-1} t H'(t), \eqno{(\ast \ast)}
$$
thus also $H' (t) = \alpha t^{-1} (F(t) - H(t))$.  Applying this technique for the c.d.f. $F_n$ of Kendall random walk  $X_n$ with the step variables $(U_k)$
$$
H(t)^n = \alpha t^{\alpha} \int_0^{1/t} s^{\alpha - 1} F_n(s) ds
$$
we see that
$$
F_n(t) =  H(t)^{n-1} \bigl[ H(t) + n\alpha^{-1} t H'(t)\bigr] = H(t)^{n-1} \bigl[ H(t) + n(F(t) - H(t))\bigr].
$$
Since $X_t = \sum_{n=0}^{\infty} X_n \mathbf{1}_{\{ N_t = n \}}$ we see that the c.d.f. $F_t$ of $X_t$ is given by
$$
F_t(x) = \sum_{n=0}^{\infty} F_n(x) \frac{(\lambda t)^n}{n!} e^{-\lambda t} =  \bigl[ 1 + \lambda t \bigl( F(x) - H(x) \bigr) \bigr] e^{- \lambda t (1 - H(x))}.
$$
Cumulative distribution function $F_{v,n}(t)$ of $v \oplus X_n$ is given in Lemma \ref{lem:2} by the formula
\begin{eqnarray*}
F_{v,n}(t) & = & \alpha^{-1} t^{1-\alpha} \frac{d}{dt} \left[ t^{\alpha} (1- u^{\alpha} t^{-\alpha})_+ H^n(t) \right] \\
 & = & \mathbf{1}_{[u,\infty)}(t) H^{n-1}(t)  \left[ H(t) + n (1- u^{\alpha} t^{-\alpha})_+ (F(t) - H(t)) \right].
\end{eqnarray*}

We need also to calculate the cumulative distribution function for the variable $u \oplus Y_n$ with the distribution $G_{u,n}$, where $Y_n$ is the Kendall random walk with the steps $(V_k)$ i.i.d. random variables with distribution $\nu$ and distribution function $G$. We see that $u \oplus Y_n$ has distribution $\delta_u \vartriangle_{\alpha} \nu^{\vartriangle_{\alpha} n}$ and the generalized characteristic function $(1- (tu)^{\alpha})_+ \Phi_{\nu}^n(t)$ for $n \geqslant 1$ and $\delta_u \vartriangle_{\alpha} \nu^{\vartriangle_{\alpha} 0} = \delta_u$ has the distribution function $G_{u,0}(t) = \mathbf{1}_{[u,\infty)}(t)$. Let $J(t) = \Phi_{\nu}(t^{-1})$, thus $J'(t) = \alpha t^{-1} (G(t) - J(t))$. Using formula ($\ast\ast$) for $\mu$ replaced by $\delta_u \diamond \nu^{\diamond n}$ we obtain for $n \geqslant 1$
\begin{eqnarray*}
G_{u,n}(t) & = & \alpha^{-1} t^{1-\alpha} \frac{d}{dt} \left[ t^{\alpha} (1- u^{\alpha} t^{-\alpha})_+ J^n(t) \right] \\
 & = & \mathbf{1}_{[u,\infty)}(t) J^{n-1}(t)  \left[ J(t) + n (1- u^{\alpha} t^{-\alpha})_+ (G(t) - J(t)) \right].
\end{eqnarray*}
This distribution has an atom at the point $u$ of the weight $G_{u,n}(u^+) = J(u)^n$ and the absolutely continuous part with the density, $t\geq u$, 
\begin{eqnarray*}
g_{u,n}(t) & = & n J^{n-1}(t)\left[ \frac{2\alpha}{t} G(t) \left(1-n\Psi\left(\frac{u}{t}\right)\right) + \Psi\left(\frac{u}{t}\right) g(t) \right]
 \\
 & - & \frac{n\alpha}{t} J^{n}(t)\left[ 2 - (n+1) \Psi\left(\frac{u}{t}\right)\right] + n(n-1) \frac{\alpha}{t} \Psi\left(\frac{u}{t}\right) J^{n-2}(t) G^2(t).
\end{eqnarray*}
\color{black}
Since $u \oplus Y_t = \sum_{n=0}^{\infty} (u \oplus Y_n) \mathbf{1}_{\{ N_t = n \}}$ we see that the c.d.f. $G_{u,t}$ of $u \oplus Y_t$ is given by
\begin{eqnarray*}
G_{u,t} (x) & = &  \sum_{n=0}^{\infty} G_{u,n} (x) \frac{(\lambda t)^n}{n!} e^{-\lambda t} \\
 & = & \mathbf{1}_{[u,\infty)}(x) \bigl[ 1 + \lambda t \bigl( 1 - u^{\alpha} x^{-\alpha} \bigr) \bigl( G(x) - J(x) \bigr) \bigr] e^{-\lambda t( 1 - J(x))}.
\end{eqnarray*}
This distribution has an atom at $u$ of the weight $G_{u,t}(u^+) = e^{-\lambda t(1 - J(u))}$.

\subsection{First safety condition for the insurance company}
In the classical theory the first safety condition for the insurance company states that $\mathbf{E} R_t >0$ for all $t > 0$. In our case we have
$$
\mathbf{E} (u \oplus Y_t)^{\alpha} - \mathbf{E} X_t^{\alpha} > 0 \quad \hbox{ for all } \quad t>0.
$$
First we calculate $\mathbf{E} X_t^{\alpha}$ assuming that the distribution of $U_1$ is absolutely continuous with respect to the Lebesgue measure (if this is not the case we shall add the atomic part):
\begin{eqnarray*}
\mathbf{E} X_t^{\alpha} & = & \int_0^{\infty} x^{\alpha} d F_t(x) = \int_0^{\infty} \alpha x^{\alpha-1}\left( 1 - F_t(x) \right) dx \\
 & = & \int_0^{\infty} \alpha x^{\alpha-1}\left[ 1 - \left( 1 + \frac{\lambda t}{\alpha} x H'(x) \right) e^{-\lambda t(1-H(x))} \right] dx \\
 & = & \int_0^{\infty} \left[ (x^{\alpha})' - \left( x^{\alpha} e^{-\lambda t(1-H(x))}\right)' \right] dx \\
 & = &  x^{\alpha} \left( 1 - e^{-\lambda t(1-H(x))}\right)\Big|_0^{\infty} = \lim_{x\rightarrow \infty} x^{\alpha} \left( 1 - e^{-\lambda t(1-H(x))}\right).
\end{eqnarray*}
In the similar way for absolutely continuous distribution of $V_1$ we obtain
\begin{eqnarray*}
\lefteqn{\mathbf{E}\bigl( u \oplus Y_t\bigr)^{\alpha} = u^{\alpha} G_{u,t}(u^+) + \int_0^{\infty} \alpha x^{\alpha-1}\left( 1 - G_{u,t}(x)\right) dx} \\
 & = & u^{\alpha}G_{u,t}(u^+) + \int_0^u \!\! \alpha x^{\alpha -1}dx + \int_u^{\infty} \!\! \left[ x^{\alpha} - \left(x^{\alpha}-u^{\alpha} \right)e^{-\lambda t (1-J(x))} \right]' dx \\
 & = &  u^{\alpha}G_{u,t}(u^+) + u^{\alpha} + \left[ x^{\alpha} - \left(x^{\alpha} -u^{\alpha} \right)e^{-\lambda t (1-J(x))} \right] \big|_u^{\infty} \\
 & = & u^{\alpha}e^{-\lambda t(1 - J(u))} +  \lim_{x\rightarrow \infty} \left[ x^{\alpha} - \left(x^{\alpha} -u^{\alpha} \right)e^{-\lambda t (1-J(x))} \right].
\end{eqnarray*}

If we consider as $\nu$ the distribution with the lack of memory property for the Kendall convolution (see \cite{Poisson}) then for some $c>0$ we have
\begin{eqnarray*}
G(x) & = & \min\{(cx)^{\alpha}, 1\}, \\
J(x) & = &  \frac{1}{2} (cx)^{\alpha} \mathbf{1}_{[0,c^{-1}]}(x) + \bigl( 1 - \frac{1}{2} (cx)^{-\alpha} \bigr) \mathbf{1}_{(c^{-1}, \infty)}(x),
\end{eqnarray*}
and, assuming that $\mathbf{1}_{[a,b]} \equiv 0$ for $a>b$ we have
\begin{eqnarray*}
G_{u,t} (x) & = & \left[ 1 + \frac{\lambda t}{2} (cx)^{\alpha} (1-u^{\alpha} x^{-\alpha} ) \right] e^{-\lambda t (1- \frac{1}{2} (cx)^{\alpha})}  \mathbf{1}_{[u,c^{-1}]}(x) \\
  && \hspace{-12mm} + \left[  1 + \frac{\lambda t}{2} (cx)^{-\alpha} (1-u^{\alpha} x^{-\alpha} ) \right] e^{- \frac{\lambda t}{2} (cx)^{-\alpha}} \mathbf{1}_{[u\vee c^{-1}, \infty)}(x).
\end{eqnarray*}
Notice that in this case
$$
G_{u,t}(u^+) = \left\{ \begin{array}{lcl}
    e^{- \lambda t (1 - \frac{1}{2} (cu)^{\alpha})} & \hbox{if} & u \leqslant c^{-1}, \\
    e^{ - \frac{\lambda t}{2} (cu)^{- \alpha}} & \hbox{if} & u > c^{-1}.
\end{array} \right.
$$
Consequently, for $u>c^{-1} = \frac{\alpha +1}{\alpha} \int x dG(x)$, which is a natural assumption since the initial capital shall be significant, we have
\begin{eqnarray*}
\mathbf{E}\bigl( u \oplus Y_t\bigr)^{\alpha}& = & u^{\alpha} e^{-\frac{\lambda t}{2} (cu)^{-\alpha} } + \lim_{x\rightarrow \infty} x^{\alpha} \left[1 -\left(1-u^{\alpha} x^{-\alpha} \right)e^{-\frac{\lambda t}{2} (cx)^{-\alpha}} \right] \\
& = &  u^{\alpha} e^{-\frac{\lambda t}{2} (cu)^{-\alpha} } +  u^{\alpha} + \frac{\lambda t}{2}\, c^{-\alpha}.
\end{eqnarray*}
For $\mu$ with the lack of memory property with $c>0$ in the Kendall convolution algebra we have 
\begin{eqnarray*}
F(x) & = &  \min\{(cx)^{\alpha}, 1\}, \\
H(x) & = &  \frac{1}{2} (cx)^{\alpha} \mathbf{1}_{[0,c^{-1}]}(x) + \bigl( 1 - \frac{1}{2} (cx)^{-\alpha} \bigr) \mathbf{1}_{(c^{-1}, \infty)}(x),
\end{eqnarray*}
Since
$$
\mathbf{E} X_t^{\alpha}  = \lim_{x\rightarrow \infty} x^{\alpha} \left( 1 - e^{-\lambda t(1-H(x))}\right) = \frac{\lambda t}{2} c ^{-\alpha}
$$
we have
$$
\mathbf{E}\bigl( u \oplus Y_t\bigr)^{\alpha} - \mathbf{E} X_t^{\alpha} = u^{\alpha} e^{-\frac{\lambda t}{2} (cu)^{-\alpha} } +  u^{\alpha} > 0,
$$
i.e. the first safety condition holds.
\subsection{Ruin probability in the infinite time horizon}

Let $Q_{\infty}(u)$ be the ruin probability for our model:
$$
Q_{\infty}(u) = 1 - \mathbf{P}\bigl\{u \oplus Y_k >  X_k  \hbox{ for all } k \in \mathbb{N} \bigr\} =: 1 - \delta(u).
$$
For the convenience we shall use the following notation: for the Markov sequence $X_n$ starting in the point $v$ we will write
$$
X^v_1 = v \oplus U_1, \,\, X_2^v = v \oplus U_1 \oplus U_2, \cdots
$$
and for the Markov sequence $Y_n$ starting at the point $u>0$ we write
$$
Y_1^u = u \oplus V_1,\,\, Y_2^u = u \oplus V_1 \oplus V_2, \cdots
$$
Let
\begin{eqnarray*}
\Lambda(v,u) & = & \mathbf{P}\bigl\{u \oplus Y_k >  v \oplus X_k  \hbox{ for all } k \in \mathbb{N} \bigr\} \\
& = & \mathbf{P}\bigl\{u \oplus V_1 > v \oplus U_1, u \oplus V_1 \oplus V_2 > v \oplus U_1 \oplus U_2 \dots \bigr\}.
\end{eqnarray*}
We need to calculate $\delta(u) = \Lambda(0,u)$. Thus
\begin{eqnarray*}
\lefteqn{\Lambda(v,u) =} \\
&& \hspace{-5mm} \int_0^{\infty} \mathbf{P}\bigl\{y_1 > v \oplus U_1, y_1 \oplus V_2 > v \oplus U_1 \oplus U_2 \dots \bigr\} \delta_u \vartriangle_{\alpha} \nu(dy_1) = \\
&&\hspace{-6mm} \int_0^{\infty} \!\! \int_{\{y_1 > x_1\}} \hspace{-5mm}\mathbf{P}\bigl\{y_1 \oplus V_2 > x_1 \oplus U_2, y_1 \oplus V_2 \oplus V_3 > x_1 \oplus U_2 \oplus U_3 \dots \bigr\} \\
&& \hspace{35mm} \delta_u \vartriangle_{\alpha} \nu(dy_1) \, \delta_v \vartriangle_{\alpha} \mu(dx_1) \\[2mm]
& = & \int_0^{\infty} \!\! \int_{\{y_1 > x_1\}} \hspace{-5mm}\Lambda(x_1,y_1) \, \delta_u \vartriangle_{\alpha} \nu(dy_1)\,  \delta_v \vartriangle_{\alpha} \mu(dx_1)
\end{eqnarray*}

Now we have for $v>u$ 
\begin{eqnarray*}
\lefteqn{\Lambda (v,u) = } \\
&&\hspace{-5mm} \int_{\lfloor v}^{\infty} \!\! \int_{x}^{\infty} \hspace{-4mm}\Lambda(x,y) \, d \left[   \Psi \left( \frac{u}{y} \right) G(y) + \left(1-\Psi \left( \frac{u}{y} \right)\right) J(y)\right]  \delta_v \vartriangle_{\alpha} \mu(dx),
\end{eqnarray*}
and if $v\leqslant u$
\begin{eqnarray*}
\lefteqn{\Lambda (v,u) =  J(u) \int_{\lfloor v}^{\infty} \Lambda(x,u) \, \delta_v \vartriangle_{\alpha} \mu(dx) } \\
  &&\hspace{-10mm} + \int_{\lfloor v}^{\infty}  \int_{u}^{\infty} \hspace{-4mm}\Lambda(x,y) \, d \left[   \Psi \left( \frac{u}{y} \right)\! G(y) + \!\left(1-\Psi \left(\frac{u}{y} \right)\!\right) J(y)\right] \delta_v \vartriangle_{\alpha} \mu(dx).
\end{eqnarray*}
This final formula shall be treated as the integral functional equation for the function $\Lambda$ and its solution (in general or depending on the specified step distributions $\mu$ and $\nu$) is an open question. 
\\
{\bf Acknowledgements.} This paper is a part of project "First order Kendall maximal autoregressive processes and their applications", \\ Grant no POIR.04.04.00-00-1D5E/16, which is carried out within the POWROTY/REINTEGRATION programme of the Foundation for Polish Science co-financed by the European Union under the European Regional Development Fund.

\end{document}